\documentclass[12pt]{article}  
\usepackage{amsthm,amsfonts,amsmath,amscd,amssymb,latexsym}   
\usepackage{epsfig,graphicx,color}    
\usepackage{color}   
\usepackage[all]{xy} 
\usepackage[titles]{tocloft} 
\title{Homogeneous quasimorphisms \\
on the symplectic linear group}  

\author{Gabi Ben Simon \qquad Dietmar~A.~Salamon \\
        ETH-Z\"urich}

\date{11 December 2007}  
 
\newcommand{\R}{{\mathbb{R}}}   
\newcommand{\cP}{{\mathcal{P}}}  
\newcommand{\Sp}{{\mathrm{Sp}}}   
\newcommand{\tSp}{{\widetilde{\Sp}}} 
\newcommand{\rG}{{\mathrm{G}}}   
   
\newcommand{\rO}{{\mathrm{O}}}   
\newcommand{\rU}{{\mathrm{U}}}   
\newcommand{\one}   
{{{\mathchoice \mathrm{ 1\mskip-4mu l} \mathrm{ 1\mskip-4mu l}   
\mathrm{ 1\mskip-4.5mu l} \mathrm{ 1\mskip-5mu l}}}}   
\def\Abs#1{\left|#1\right|}   
   
\begin{document}  

\maketitle   
 
%%%%%%%%%%%%%%%%%%%%%%%%%%%%%%%%%%%%%%%%%%%%%%   
%%%%%%%%%%%%%%%%%%%%%%%%%%%%%%%%%%%%%%%%%%%%%%   
%%%%%%%%%%%%%%%% Abstract %%%%%%%%%%%%%%%%%%%%   
%%%%%%%%%%%%%%%%%%%%%%%%%%%%%%%%%%%%%%%%%%%%%%   
%%%%%%%%%%%%%%%%%%%%%%%%%%%%%%%%%%%%%%%%%%%%%%   
   
%begin{abstract}   
%end{abstract}   

Let $\rG$ be a group. A {\bf quasimorphism} on $\rG$ 
is a map $\rho:\rG\to\R$ satisfying 
$$
\Abs{\rho(gh)-\rho(g)-\rho(h)}\le C
$$ 
for all $g,h\in\rG$ and a suitable constant~$C$. 
It is called {\bf homogeneous} if $\rho(g^k)=k\rho(g)$ 
for every $g\in\rG$ and every integer $k\ge 0$. 
Let
$$
\Sp(2n) := \left\{\Psi\in\R^{2n\times 2n}\,|\,
\Psi J_0\Psi^T = J_0\right\},\qquad
J_0 := \left(\begin{array}{rr}
0 & -\one \\ \one & 0
\end{array}\right),
$$ 
denote the group of symplectic matrices 
and $\tSp(2n)$ its universal cover. Think of an element 
of $\tSp(2n)$ as a homotopy class $[\Psi]$ (with fixed endpoints)
of a smooth path $\Psi:[0,1]\to\Sp(2n)$ satisfying $\Psi(0)=\one$.

\medskip\noindent
{\bf Theorem~1.}
{\it There is a unique homogeneous quasimorphism $\mu$ on $\tSp(2n)$
that descends to the determinant homomorphism on $\rU(n)$ 
in the sense that
$$
\det(X+iY) = \exp\left(2\pi i\mu([\Psi])\right),\qquad
\left(\begin{array}{rr}
X & -Y \\ Y & X
\end{array}\right) := \Psi(1),
$$
for every $[\Psi]\in\tSp(2n)$ with 
$\Psi(1)\in\Sp(2n)\cap\rO(2n)\cong\rU(n)$.}

\medskip\noindent 
The quasimomorphism of Theorem~1 plays a central role 
in~\cite{BS} and this motivated the present note. 
Two explicit constructions of the quasimorphism 
can be found in~\cite{BG} and~\cite{SZ}.  The construction in~\cite{BG} 
uses the unitary part in a polar decomposition and homogenization.
The construction in~\cite{SZ} uses the eigenvalue decomposition 
of a symplectic matrix (but does not mention the 
term {\it quasimorphism}). 

\medskip\noindent
{\bf Lemma~1}
{\it If $\rho:\rG\to\R$ is a homogeneous quasimorphism
then $\rho$ is invariant under conjugation and 
$\rho(g^{-1}) = -\rho(g)$ for every $g\in\rG$.}

\begin{proof}[Proof of Lemma~1]
Let $C$ be the constant in the definition of quasimorphism.
By homgeneity, we have $\rho(1)=0$. Hence 
$\Abs{\rho(g^k)+\rho(g^{-k})}\le C$ for every ${g\in\rG}$
and every integer $k\ge 0$. By homogeneity,
we obtain ${\Abs{\rho(g)+\rho(g^{-1})}\le C/k}$ for 
every $k$ and so $\rho(g^{-1}) = -\rho(g)$.
Hence 
$$
\Abs{\rho(ghg^{-1})-\rho(h)}
= \Abs{\rho(ghg^{-1})-\rho(g)-\rho(h)-\rho(g^{-1})}
\le 2C.
$$
Using homogeneity again we obtain $\rho(ghg^{-1})=\rho(h)$
for all $g,h\in\rG$.  
\end{proof}

\begin{proof}[Proof of Theorem~1]
Let $\cP\subset\Sp(2n)$ denote the set of 
symmetric positive definite symplectic matrices. This space 
is contractible and hence there is a natural injection
$\iota:\cP\to\tSp(2n)$. Explicitly, the map $\iota$ assigns 
to a matrix $P\in\cP$ the unique homotopy class of paths
$\Phi:[0,1]\to\cP$ with endpoints $\Phi(0)=\one$ 
and $\Phi(1)=P$.  

Let $\mu:\tSp(2n)\to\R$ be a homogeneous 
quasimorphism that descends to the determinant 
homomorphism on $\rU(n)$. It suffices to prove that 
the restriction of $\mu$ to $\iota(\cP)$ is bounded.
(If $\mu'$ is another quasimorphism satisfying the 
requirements of Theorem~1 and $\mu,\mu'$ 
are bounded on $\iota(\cP)$ then, 
by polar decomposition and the determinant assumption,
their difference is bounded and so, by homogeneity, they are equal.) 
We prove that $\mu$ vanishes on $\iota(\cP)$.
For every unitary matrix $Q\in\rU(n)\subset\Sp(2n)$ 
and every $P\in\cP$ we have 
$$
\mu(\iota(QPQ^T))=\mu(\iota(P)).\leqno{(1)}
$$
To see this, choose two paths $\Phi:[0,1]\to\cP$ and $\Psi:[0,1]\to\rU(n)$
such that $\Phi(0)=\Psi(0)=1$ and $\Phi(1)=P$, $\Psi(1)=Q$. 
Then $\mu([\Phi])=\mu([\Psi\Phi\Psi^{-1}])$, by Lemma~1,
and so~(1) follows from the fact that $\Psi^{-1}=\Psi^T$.
Now let $P\in\cP$. Since $P$ is a symmetric symplectic matrix 
we have $PJ_0P=J_0$ and hence
$$
\mu(\iota(P)) = \mu(\iota(J_0P^{-1}J_0^{-1}))
= \mu(\iota(P^{-1})) = \mu(\iota(P)^{-1}) = -\mu(\iota(P)).
$$
Here the second equation follows from~(1) and the
last from Lemma~1.  This shows that 
$\mu(\iota(P))=0$ for every $P\in\cP$.
\end{proof}

\noindent{\bf Remark~1.}
Lemma~1 is well known to the experts.  
We included a proof to give a self-contained exposition,
and because we didn't find an explicit reference.

\medskip\noindent{\bf Remark~2.}
Related results, obtained with different methods,
are contained in~\cite{BG} and~\cite{BIW}.
Our main theorem can in fact be deduced from
these results.

\medskip\noindent{\bf Remark~3.}
The determinant homomorphism $\det:\rU(n)\to S^1$
is uniquely determined by the condition that it induces 
an isomorphism on fundamental groups.  Hence it follows from 
Theorem~1 that the homogeneous quasimorphism $\mu:\tSp(2n)\to\R$
is uniquely determined by the condition that it restricts
to an isomorphism of the fundamental group of $\Sp(2n)$ 
to the integers.  

%\medskip\noindent{\bf Acknowledgement.}
%Thanks to Marc Burger and Leonid Polterovich 
%for helpful suggestions.

%%%%%%%%%%%%%%%%%%%%%%%%%%%%%%%%%%%%%%%%%%%%%%   
%%%%%%%%%%%%%%%%%%%%%%%%%%%%%%%%%%%%%%%%%%%%%%   
%%%%%%%%%%%%%%%% References %%%%%%%%%%%%%%%%%%   
%%%%%%%%%%%%%%%%%%%%%%%%%%%%%%%%%%%%%%%%%%%%%%   
%%%%%%%%%%%%%%%%%%%%%%%%%%%%%%%%%%%%%%%%%%%%%%   

\end{document}